\documentclass{agt}

\usepackage{amsfonts}
\usepackage{amsmath}
\usepackage{stmaryrd}
\usepackage{amssymb}
\usepackage{graphicx}
\usepackage{appendix}
\usepackage{subfigure}
\usepackage{url}
\usepackage{hyperref}

\begin{document}
\title{Computational Topology Counterexamples with \\ 3D Visualization of B\'ezier Curves}{Computational Topology Counterexamples}
\author{J. Li \and T. J. Peters\thanks{This author was partially supported by NSF grants CCF 0429477, CMMI 1053077 and CNS 0923158, as well as by an IBM Faculty Award and IBM Doctoral Fellowships.  All statements here are the responsibility of the author, not of the National Science Foundation nor of IBM.} \and D. Marsh \and K. E. Jordan}{}

\address[ji.li@uconn.edu]{J. Li}{Department of Mathematics, University of Connecticut, Storrs, CT, USA.}
\address[tpeters@cse.uconn.edu]{T. J. Peters}{Department of Computer Science and Engineering,
        University of Connecticut, Storrs, CT, USA.}
\address[david.the.marsh@gmail.com]{D. Marsh}{Pratt and Whitney, East Hartford, CT, USA.}
\address[kjordan@us.ibm.com]{K. E. Jordan}{IBM T.J. Watson Research, Cambridge Research Center, Cambridge, MA, USA.}

\date{\today}

\begin{abstract} 
For applications in computing, B\'ezier curves are pervasive and are defined by a piecewise linear curve $\mathcal{L}$ which is embedded in $\mathbb{R}^3$ and yields a smooth polynomial curve $\mathcal{C}$ embedded in $\mathbb{R}^3$.  It is of interest to understand when $\mathcal{L}$ and $\mathcal{C}$ have the same embeddings.  One class of counterexamples is shown for $\mathcal{L}$ being unknotted, while $\mathcal{C}$ is knotted.   Another class of counterexamples is created where $\mathcal{L}$ is equilateral and simple, while $\mathcal{C}$ is self-intersecting.   These counterexamples were discovered using curve visualizing software and numerical algorithms that produce general procedures to create more examples. 
\end{abstract}

\begin{figure}[h!]
\centering
        \subfigure[Unknotted $\mathcal{L}$ with knotted $\mathcal{C}$]
   {   \includegraphics[height=3.5cm]{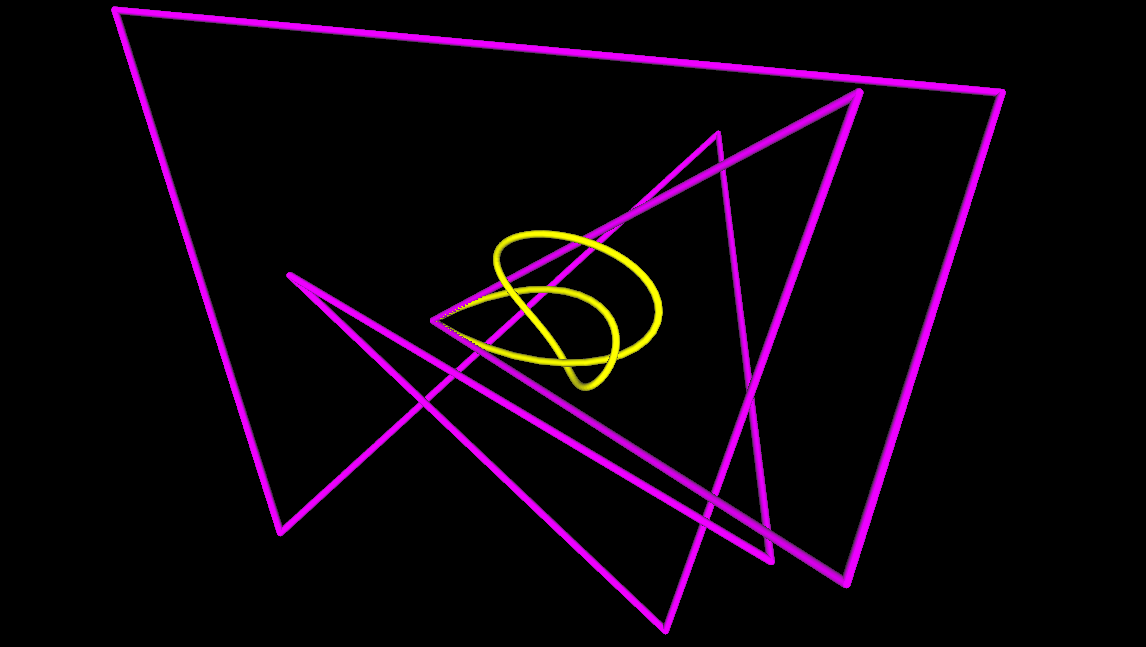}  \label{fig:k} }
                         \subfigure[Zoomed-in view of $\mathcal{C}$]
   {   \includegraphics[height=3.5cm]{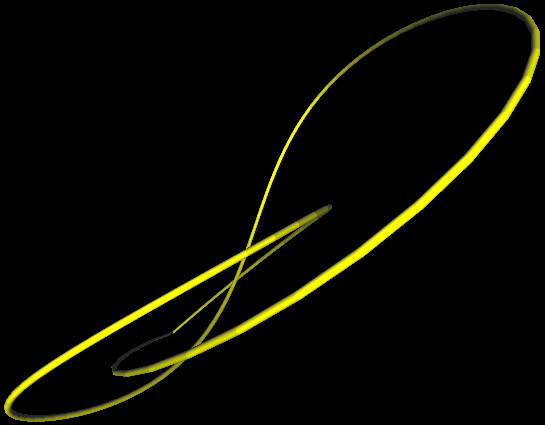}  \label{fig:kz} }
    \caption{Unknotted $\mathcal{L}$ with knotted $\mathcal{C}$}    \label{fig:knots}    
 \end{figure}   
   
\begin{figure}[h!]   
       \subfigure[ Equilateral, simple $\mathcal{L}$ with self-intersecting $\mathcal{C}$]
    {   \includegraphics[height=3.5cm]{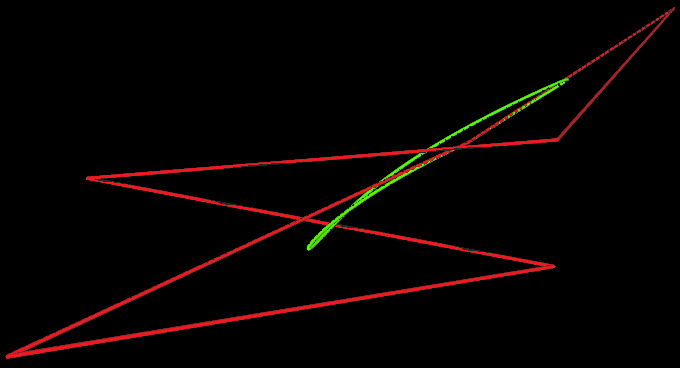}   \label{fig:eq}   }
               \subfigure[Zoomed-in self-intersection]
   {   \includegraphics[height=3.5cm]{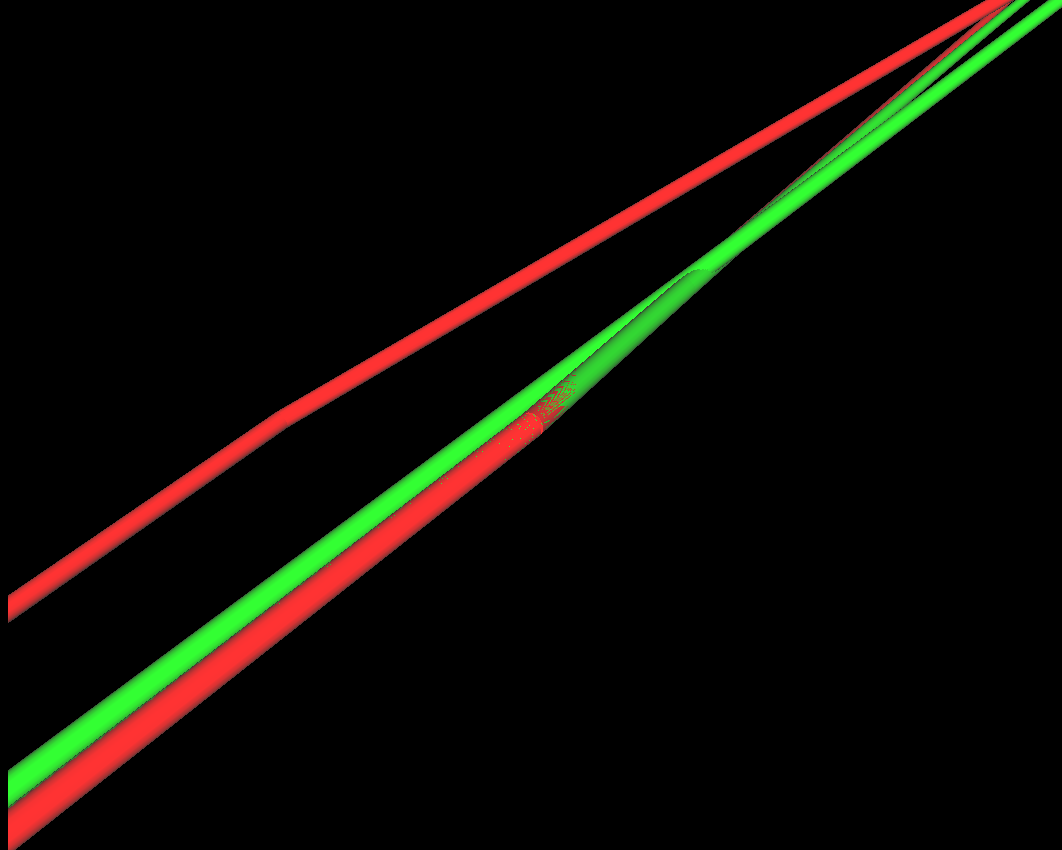}  \label{fig:eqz} }

    \caption{Equilateral, simple $\mathcal{L}$ with self-intersecting $\mathcal{C}$}    \label{fig:knots1}    
 \end{figure}

\section{Introduction}
\label{sec:intro}

Computational topology lies properly within the broad scope of applied general topology and depends upon a novel integration of the `pure mathematics' of topology with the `applied mathematics' of numerical analysis.  Computational topology also blends general topology, geometric topology and knot theory with computer visualization and graphics, as presented here.   For the computational cases considered here, ideas from general topology, geometric topology and knot theory are complemented by numeric arguments in a novel integration of the `pure' field of topology with the `applied' focus of numerical analysis.  Additionally, aspects of computer visualization and graphics are incorporated into the proofs.  Much of this work was motivated by modeling biological molecules, such as proteins, as knots for visualization synchronized to simulations of the writhing of these molecules.

Attention is restricted to when $\mathcal{C}$ is closed, implying that $\mathcal{L}$ is also closed.  As $\mathcal{C}$ is created from $\mathcal{L}$, it is natural to ask which topological characteristics are shared by these two curves, particularly as the control polygon often serves as an approximation to the B\'ezier curve in many practical applications. However, topological differences between a B\'ezier curve and its control polygon can exist and it is natural to develop counterexamples to show these topological differences. These counterexamples extend beyond related results, while we introduce new computational methods to generate additional counterexamples.

In Section~\ref{sec:ex1}, we present a counterexample of B\'ezier curve and its control polygon being homeomorphic, yet  not ambient isotopic. To develop this counterexample, we created and used a computer visualization tool to study topological relationships between a B\'ezier curve and its control polygon.  We viewed the images to motivate formal proofs, which also rely on numerical analysis and geometric arguments.

In Section~\ref{sec:eqsimnot}, we present numerical techniques to create a class of topological counterexamples -- where a B\'ezier curve and its control polygon are not even homeomorphic, as the B\'ezier curve is self-intersecting while the control polygon is simple.  We exhibit self-intersection by a numerical method, which finds the roots of a system of equation.  We freely admit that these roots are not determined with infinite precision, but such calculations on polynomials of degree 6, as in these examples, typically elude precise calculation.  We argue that two primary values of our method for these approximated solutions are
\begin{enumerate}
\item as a catalyst to alternative examples that may admit infinite precision calculation and rigorous topological proofs, and
\item as having specified digits of accuracy -- typically crucial for acceptable approximations in computational mathematics and computer graphics.
\end{enumerate} 

Using the visualization tool described in Section~\ref{sec:vizb}, we viewed many examples where the B\'ezier curve was simple, while its control polygon was equilateral and simple.  It is well known that B\'ezier curve can be self-intersecting even when its control polygon is simple, but we conjectured that the added equilateral hypothesis would imply that both curves were simple.  
While this visual evidence was suggestive, we present a general numerical approach in Section~\ref{sec:eqsimnot} that supports a contrary interpretation.  This example provides guidance for designing appropriate approximation algorithms for computer graphics.

\subsection{Mathematical Definitions}
\label{ssec:mdef}

The definitions presented are restricted to $\mathbb{R}^3$, as sufficient for the purposes of this paper, but the interested reader can find appropriate generalizations in published literature.

\begin{defn}
\label{def:knot}
A knot \cite{Armstrong1983} is a curve in $\mathbb{R}^3$ which is homeomorphic to a circle.
\end{defn}

Knots are often described by a {\em knot diagram} \cite{Livingston1993}, which is a planar projection of a knot.  Self-intersections in the knot diagram correspond to {\em crossings} in the knot, where each crossing has one arc that is an {\em undercrossing} and an {\em overcrossing}, relative to the direction of projection.

Homeomorphism is an equivalence relation over point sets and does not distinguish between different embeddings, while ambient isotopy is a stronger equivalence relation which is fundamental for classification of knots.

\begin{defn}
Let $X$ and $Y$ be two subspaces of $\mathbb{R}^3$.  A continuous function
\[ H:\mathbb{R}^3 \times [0,1] \to \mathbb{R}^3 \]
is an {\bf ambient isotopy} between $X$ and $Y$ if $H$ satisfies the following:
\begin{enumerate}

\item $H(\cdot, 0)$ is the identity,

\item $H(X,1) = Y$, and

\item $\forall t \in [0,1], H(\cdot,t)$ is a homeomorphism from
$\mathbb{R}^3$ onto $\mathbb{R}^3$.

\end{enumerate}
The sets $X$ and $Y$ are then said to be {\bf ambient isotopic}.
\label{def:aiso}
\end{defn}

\begin{defn}\label{def:c}
Denote $\mathcal{C}(t)$ as the parameterized B\'ezier curve of degree $n$ with control points $P_m \in \mathbb{R}^3$, defined by
$$
\mathcal{C}(t)=\sum_{i=0}^{n}{B_{i,n}(t)P_i},
t\in[0,1]
$$
where $B_{i,n}(t) = \left(\!\!\!
  \begin{array}{c}
	n \\
	i
  \end{array}
  \!\!\!\right)t^i(1-t)^{n-i}$.
\end{defn}

\subsection{Related Literature}
\label{ssec:relwork}

A B\'ezier curve and its control polygon may have substantial topological differences. It is well known that a B\'ezier curve and its control polygon are not necessarily homeomorphic \cite{Piegl}.  Recently, it was shown that there exists an unknotted B\'ezier curve with a knotted control polygon  \cite{Bisceglio}.  A specific dual example has also been shown \cite{Carlo} of a knotted B\'ezier curve with an unknotted control polygon.  However, the methodology was a visual construction without formal proof and is not easily generalized.  Software, {\em Knot\_Spline\_Vis}, developed by authors Marsh and Peters, was used to find another example, where the methodology can more easily be generalized and {\em Knot\_Spline\_Vis} is publicly available \cite{tjp-web}.    

Much of the motivation for considering these counterexamples came from applications in scientific visualization \cite{JKSP10,TCS08,JMPR11}.  A primary focus was to establish tubular neighborhoods of knotted curves so that piecewise linear (PL) approximations of those curves within those neighborhoods remained ambient isotopic to the original curves.  This was initially considered for approximations used in producing {\em static} images of these curves, but it became readily apparent that these same neighborhoods also provided bounds within which many perturbations of these models remained ambient isotopic under these movements.  That theory is being applied to dynamic visualizations of molecular simulations, where the neighborhood boundaries permit convenient warnings as to an impending self-intersection, as of possible interest to biologists and chemists who are running the simulations.

Previous work \cite{KnotPlot} in knot visualization provides a rich set of data for PL knots, where each edge is of the same length.  The interface of 
{\em Knot\_Spline\_Vis} was designed to import this equilateral PL knot data and then generate the associated B\'ezier curves.  This matured into an empirical study of dozens of cases, where all examples examined yielded simple B\'ezier curves for these simple equilateral control polygons.  This raised the question of whether the presence of equilateral edges in the control polygon might be a sufficient additional hypothesis to ensure homeomorphic equivalence with the B\'ezier curve, as the previously cited examples \cite{Piegl} did not have equilateral control polygons.  

\section{Unknotted $\mathcal{L}$ with knotted $\mathcal{C}$}
\label{sec:ex1}

\begin{figure}[h!]
\centering
    \subfigure[Unknotted $\mathcal{C}$ \& The $\mathcal{P}$]
    {
   \includegraphics[height=5cm]{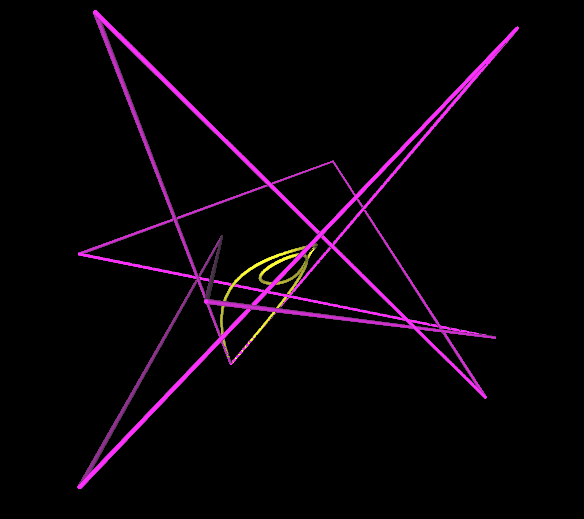}
   \label{fig:unkCP}
    }
    \subfigure[Zoomed-in View of $\mathcal{C}$]
    {
   \includegraphics[height=5cm]{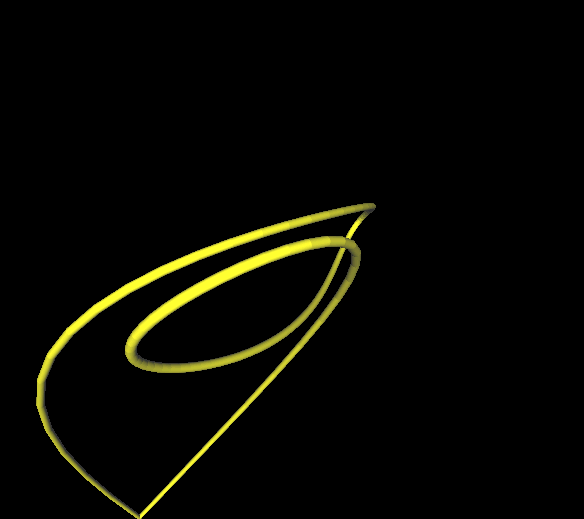}
   \label{fig:unkview}
   }
    \caption{Visual experiments}
    \label{fig:pre-unknotpoly}
\end{figure}

In order to produce a knotted B\'ezier curve with an unknotted control polygon, we invoked  {\em Knot\_Spline\_Vis} with an example (Figure~\ref{fig:pre-unknotpoly}) of an unknotted B\'ezier curve, where the total curvature appeared to be larger than $4\pi$. (The total curvature being larger than $4\pi$ is a necessary condition of knottedness.) 
We experimented on this example by moving control points to construct a B\'ezier curve that visually appeared to be knotted. We then moved control points to unknot the control polygon while keeping the B\'ezier curve knotted. In the end we obtained a B\'ezier curve and the control polygon (of degree $10$) (Figure~\ref{fig:k}) defined by the control points $\{P_0,P_1,\cdots,P_9,P_0\}$ listed below:
{\footnotesize $$(-5.9,4.7,-6.2), (10.3,-1.1,8.9), (-2.6,-12.4,-6.3), (-10,7,-0.3),$$
$$ (1.9,-12,-0.6), (11.2,7.5,-7.6), (-15.3,-1.7,-4.1), (-11.7,20,3.5), $$
$$(17.9,-1.1,2.9), (2.9,-13.7,4.8), (-5.9,4.7,-6.2).$$}

The 3D visualization offers only suggestive evidence that the above B\'ezier curve is knotted while the control polygon is unknotted. We provide rigorous mathematical proofs of these properties in Sections~\ref{ssec:p1a} and \ref{ssec:p1b}. Generally, proving knottedness or unknottedness can be difficult \cite{Hass1999}, but is accessible for the counterexample here.

\subsection{Proofs of the B\'ezier curve being knotted}
\label{ssec:p1a}

We prove that the B\'ezier curve is a trefoil\footnote{A trefoil is a knot with three alternating crossings\cite{Livingston1993}.}. We orthogonally project the curve onto $x$-$y$ plane. We then shown that there are three self-intersections in the projection and these intersections are alternating crossings in 3D.  Since projections preserve self-intersections, this curve can have no more than 3 self-intersections, but these self-intersections in $x$-$y$ plane are shown to have pre-images that are 3D crossings, so the original curve has no self-intersections. 

\begin{figure}[h!]
\begin{center}
\includegraphics[height=7cm]{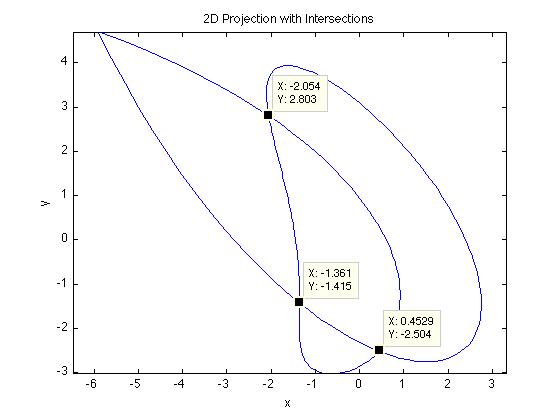} 
\end{center}

\caption{The 2D projection of the knot}\label{fig:k-2d}
\end{figure}

Since the 2D curve in Figure~\ref{fig:k-2d} has degree $10$, we use a numerical method implemented by MATLAB function `fminsearch' to find the parameters where the curve intersects with itself.  We provide the numerical codes in Appendix~\ref{app:selfint}, and we provide the data used to find these parameters in Appendix~\ref{data:selfint}. The pairs of parameters (labeled in order as $t_1,\ldots,t_6$) of the self-intersections are listed below:

$$[t_1=0.0306,t_4=0.5573], [t_2=0.1573,t_5=0.9244], [t_3=0.3731,t_6=0.9493].$$

Next we prove that these 2D intersections are projections from three alternating crossings in 3D.  The above parameters are substituted into the B\'ezier curve (Definition~\ref{def:c}) to get pairs of points (numerical codes for this calculation are in Appendix~\ref{app:cr}):
{\small $$[\mathcal{C}(t_1)=(-2.0539,2.8001,-2.6929), \mathcal{C}(t_4)=(-2.0530,2.7987,-2.0143)],$$ 
$$[\mathcal{C}(t_2)=(0.4376,-2.5212,-0.0576),\mathcal{C}(t_5)=(-0.4364,-2.5206,-0.5547)],$$ 
$$[\mathcal{C}(t_3)=(-1.3613,-1.4239,-2.2944),\mathcal{C}(t_6)=(-1.3624,-1.4232,-1.9067)].$$}

The alternating crossings follow from comparing the z-coordinates in each pair. Precisely, according to the parameters given above, the tracing of the six points in order is 
$$\mathcal{C}(t_1), \mathcal{C}(t_2), \mathcal{C}(t_3), \mathcal{C}(t_4), \mathcal{C}(t_5), \mathcal{C}(t_6),$$
and the crossings at these points are
$$under, over, under, over, under, over.$$

\subsection{Proof of the control polygon being unknotted}
\label{ssec:p1b}

To prove that the control polygon $P=(P_0,P_1,\cdots,P_9,P_0)$ is an unknot, it is necessary to show that $P$ is simple and unknotted.  We directly tested each pair of segments of $P$ for  non-self-intersection and those calculations can be repeated by the interested reader.   We prove unknottedness using a 3D {\em push}, as the obvious generalization of a 2D function from low-dimensional geometric topology \cite{Bing1983}.    We restrict attention to a {\em median push}, as defined below.   The full sequence of 5 median pushes is explicated, where the first 4 median pushes are equivalently described by Reidemeister moves \cite{Livingston1993}.

\begin{figure}[h!]
\centering
    \subfigure[The initial control polygon]
    {
   \includegraphics[height=4cm]{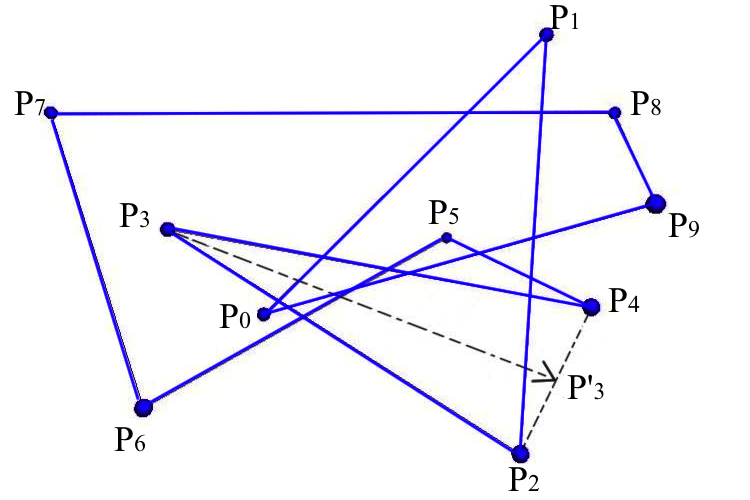}
    \label{fig:unknot0}
    }
    \subfigure[After the 1st push]
    {
   \includegraphics[height=4cm]{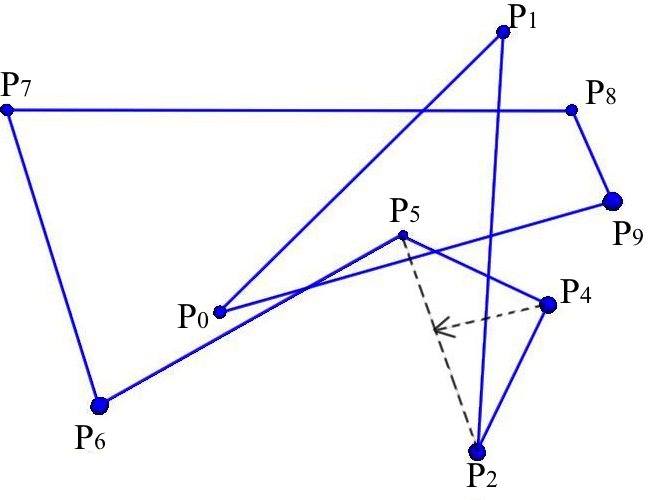}
    \label{fig:unknot1}
    }


    \subfigure[After the 2nd push]
    {
   \includegraphics[height=4cm]{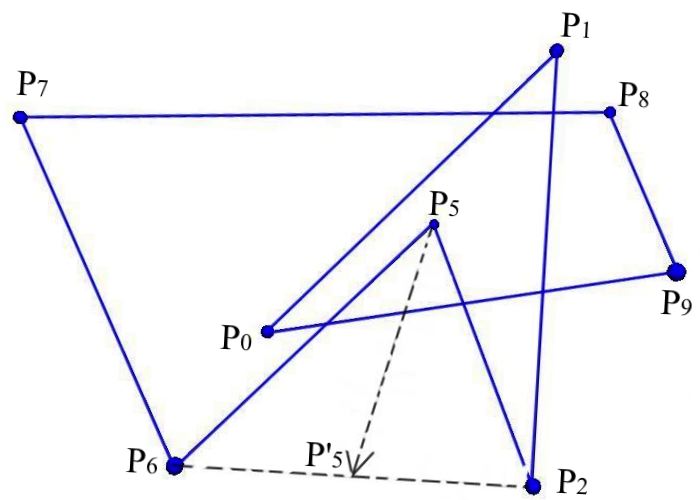}
    \label{fig:unknot2}
    }  
    \subfigure[After the 3rd push]
    {
   \includegraphics[height=4cm]{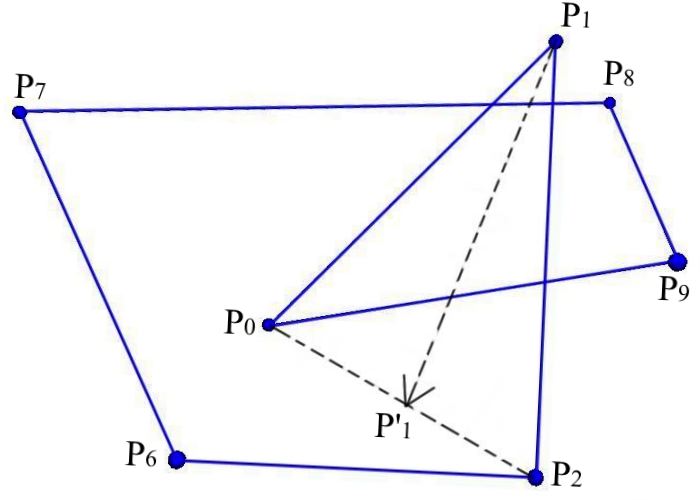}
    \label{fig:unknot3}
    }


    \subfigure[After the 4th push]
    {
   \includegraphics[height=4cm]{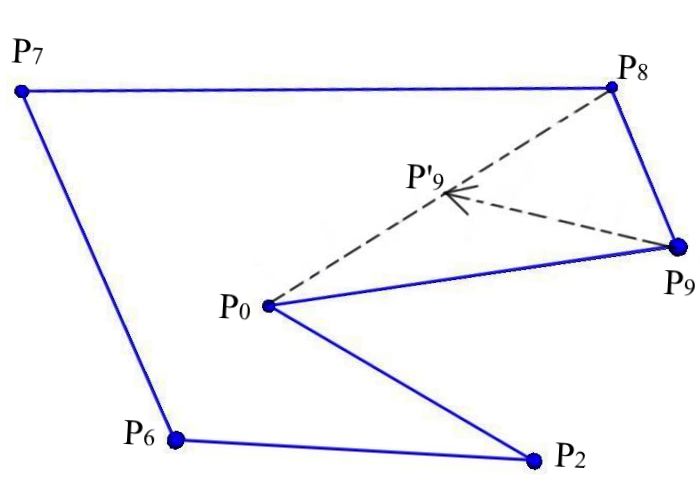}
    \label{fig:unknot4}
    }
    \subfigure[The unknot after these pushes]
    {
   \includegraphics[height=4cm]{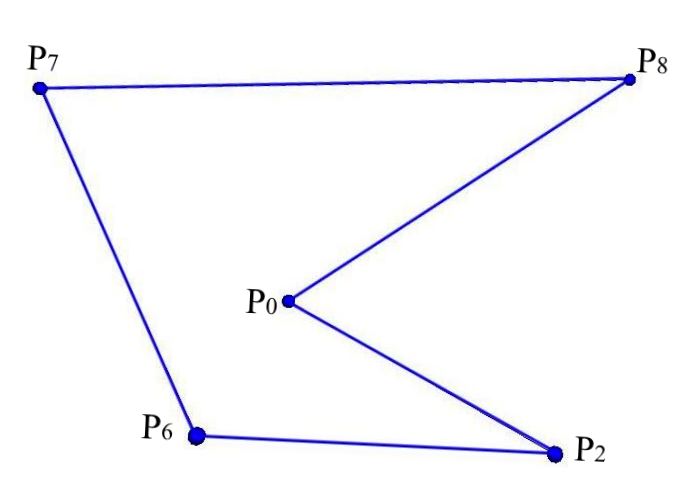}
    \label{fig:unknot5}
    }
     \caption{Pushes and PL curves in 3D}
    \label{fig:pushes}
\end{figure}

\begin{defn}
\label{def:mp}
Suppose $\triangle{P_{j-1}P_jP_{j+1}}$ is a triangle determined by non-collinear vertices $P_{j-1}, P_j$ and $P_{j+1}$ of $P$. Push the vertex $P_j$, along the corresponding median of the triangle to the middle point of the side $P_{j-1}P_{j+1}$. We call this function a median push.\end{defn}

We depict the sequence of pushes used in 
Figure~\ref{fig:pushes}, showing the 3D graphs of the PL space curves after the median pushes. These graphs show, at each step, which vertex is pushed and its image. For example, in Figure~\ref{fig:unknot0}, the vertex $P_3$ is pushed to $P'_3$ and the resultant polygon after this push is shown by Figure~\ref{fig:unknot1}.  Figures~~\ref{fig:unknot0}~\ref{fig:unknot1}~\ref{fig:unknot2} and~\ref{fig:unknot3} have corresponding Reidemeister moves, as those pushes eliminate at least one crossing.  Using the published notation \cite{Livingston1993}, Figure~\ref{fig:unknot0} depicts a Reidemeister move of Type 2b.  Similarly, Figure~\ref{fig:unknot1} depicts a Reidemeister move of Type 1b; Figure~\ref{fig:unknot2} has Type 2b and Figure~\ref{fig:unknot3} has a move of Type 2b, followed by a move of Type 1b.  The final push to achieve Figure~\ref{fig:unknot5} does not correspond to any Reidemeister move as no crossings are changed, but it is included to have a polygon with only five edges, which necessarily must be the unknot \cite{adams2004knot}.

We prove that $P$ is unknotted by showing that $P$ is ambient isotopic to the unknotted PL curve shown in Figure~\ref{fig:unknot5}. As a sufficient condition, we show that the pushes do not cause intersections \cite{J.W.Alexander_G.B.Briggs1926}. (We gained significant intuition for specifying the sequence of pushes by visual verification with our 3D graphics software capabilities to translate, rotate and zoom the images.)   We now present the formal arguments.

Consider Figure~\ref{fig:unknot0}, where $P_3$ is pushed to $P'_3$ (the middle point\footnote{It is not necessary to push $P_3$ to the middle point. Any point along $\overrightarrow{P_2P_4}$ would suffice.} of $\overrightarrow{P_2P_4}$).  We show that any segments other than $\overrightarrow{P_2P_3}$ and $\overrightarrow{P_3P_4}$ of $P$ do not intersect the triangle $\triangle{P_2P_3P_4}$ or the triangle interior.

We parameterize each segment by:
$$\overrightarrow{P_iP_{i+1}}: P_i+(P_{i+1}-P_i)t, \ t\in[0,1]$$ for $i=0,1,\cdots,9$ and let $P_{10}=P_0$.   Then the points given by 
$$P_3+a(P_2-P_3)+b(P_4-P_3),$$ 
for $a,b \geq 0$ and $a+b \leq 1$ are on the $\triangle{P_2P_3P_4}$ and  contained in its interior.   Hence $P_i+(P_{i+1}-P_i)t$ intersects $\triangle{P_2P_3P_4}$ or its interior if and only if
\begin{align}\label{eq:triangle}P_i+(P_{i+1}-P_i)t=P_3+a(P_2-P_3)+b(P_4-P_3)\end{align}
has a solution for some $t\in [0,1]$ and $a,b \geq 0$ and $a+b \leq 1$.

For each $i=0,1,\cdots,9$, we solve Equation~\ref{eq:triangle} (a system of $3$ linear equations) for $a,b$ and $t$ with the above constraints. Calculations show there is no solution for each system, and hence $P_i+(P_{i+1}-P_i)t$ does not intersect $\triangle{P_2P_3P_4}$ or its interior for each $i=0,1,\cdots,9$. Thus it follows that the push does not cause any intersections. 

Similar computations verify that the other pushes do not cause any intersections, thus establishing the ambient isotopy.

\section{Equilateral, simple $\mathcal{L}$ with self-intersecting $\mathcal{C}$}
\label{sec:eqsimnot}

\begin{figure}[h!]
\begin{center}\includegraphics[height=7cm]{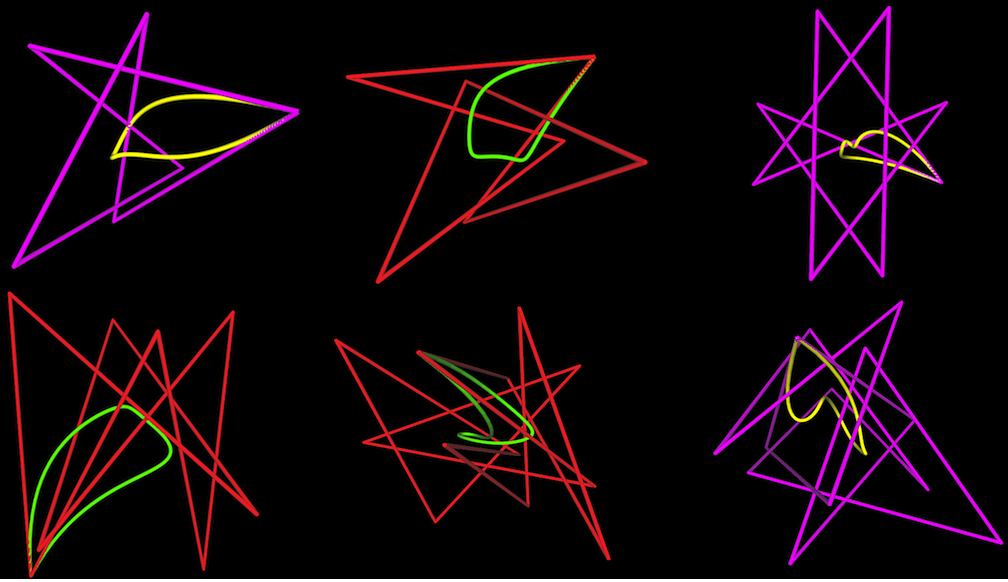}\end{center}

\caption{Unknotted B\'ezier curves with equilateral control polygons}\label{fig:examples-equilateral}
\end{figure}

We visualized many cases of simple, closed equilateral control polygons\footnote{Throughout this section, all control polygons are simple.} that all appeared to have unknotted B\'ezier curves.  Many of these control polygons were nontrivial knots.   Prompted by this visual evidence, we conjectured that any simple, closed equilateral control produces a simple B\'ezier curve, where some examples are shown in Figure~\ref{fig:examples-equilateral}.

We  now present numerical evidence to the contrary.  As noted in Section~\ref{sec:intro}, the degree 6 polynomials make precise computation difficult, so we do not provide a completely formal proof, independent of numerical methods.  A legitimate concern is whether the numerical approximation produces coordinates for a `near' intersection, subject to the accuracy of the floating point computation.
We cannot refute that possibility, but we argue that in the context of graphical images, the level of approximation produced is often sufficient.  In particular, we have parameters in the code to adjust the number of digits of accuracy.  This is level of user-defined precision is often accepted as sufficient for visualization \cite{JKSP10}.  The user can then set the graphical resolution so that points that are determined to be the same within some acceptable numerical tolerance will also appear within the same pixel.

\subsection{Intuitive overview}
\label{ssec:iover}

We create many examples to test and retain only those that satisfy all three criteria listed in Section~\ref{ssec:num}.
We begin the creation of a closed equilateral polygon by setting $P_0=(0,0,0)$. We then take $\{q_0,q_1,\cdots,q_{n-1}\}$ (Equation~\ref{eq:q}) from the unit sphere so that 
$$P_1=P_0+q_0,\ \ \ P_2=P_1+q_1,\ \ \  \cdots,\ \ \  P_n=P_{n-1}+q_{n-1}.$$
We ensure that the polygon is closed in Section~\ref{ssec:num} item $(3)$. 

We consider a sufficient and necessary condition for a B\'ezier curve being self-intersecting (Equation~\ref{eq:S2}). Since we want the equilateral polygon to define a B\'ezier curve that is self-intersecting, we not only select $\{q_0,q_1,\cdots,q_{n-1}\}$ from the unit sphere as above, but also select them such that Equation~\ref{eq:S2} is zero for some parameters $s$ and $t$. Consequently the set of control points generated determines a closed equilateral control polygon and a self-intersecting B\'ezier curve. 

\subsection{Necessary and sufficient condition for self-intersection}

We rely upon the following published equation \cite{Andersson1998} for necessary and sufficient conditions for self-intersection of a B\'ezier curve
\begin{equation}
S(s,t)=\frac{1}{n}\frac{\mathcal{C}(1-s)-\mathcal{C}(t)}{(1-s)-t},
\label{eq:Bselfint}
\end{equation}
with the domain $D=\{(s,t): s+t < 1, s,t \geq 0, (s,t) \neq (0,0)\}$. A B\'ezier curve defined by $\mathcal{C}(t)$ is self-intersecting if and only if there exist $s$ and $t$ in the domain $D$ such that $S(s,t)=0$, , with an alternative formulation\footnote{One needs to write out the \cite[Equation (6)]{Andersson1998} to obtain Equation~\ref{eq:S2}.} given by 

\begin{align}\label{eq:S2} S(s,t) = \frac{1}{n} \sum_{i=0}^{n-1} \sum_{j=0}^{n-1-i} \binom{n-1-i}{j} s^{n-1-i-j} (1-s)^j \sum_{k=0}^i \binom{i}{k}(1-t)^{i-k}t^kq_{j+k},\end{align}
where \begin{align}\label{eq:q}q_i=P_{i+1}-P_{i}\ \ \ \ for\ i=\{0,1,\cdots,n-1\}.\end{align}

\subsection{A representative counterexample generated}
\label{sec:ex}

We present a single numerical counterexample, noting that while only one counterexample is presented, the numerical algorithm implemented can be used to find many such examples.  We list six distinct control points (the seventh control point is equal to the first control point) that determine an equilateral simple $\mathcal{L}$ and a self-intersecting $\mathcal{C}$ (shown in Firgure~\ref{fig:knots1}), as generated by the algorithm described in detail in Section~\ref{ssec:num}. 

$$(0,0,0),  (0.0305,0.0810,0.9962),  (-0.2074, -0.2671,1.9030), $$
$$ (-0.1792,-0.3402,0.9063),  (0.0189, 0.0782,0.0185), (0.1557, 0.2329,-0.9600) .$$

We verify that $\mathcal{L}$ is simple by considering all pairwise intersections of  the segments of this control polygon.
The self-intersection of the B\'ezier curve occurs at a point that is numerically approximated as
$$[s,t]=[0.2969,0.0633]$$
where correspondingly,
$$S(s,t)=(-0.0003861,-0.000097,0.0001462) \approx (0,0,0).$$
The error occurs because of numerical round off on $s,t$ and the control points.

\subsection{The numerical method for generating counterexamples}
\label{ssec:num}

We provide the numerical codes and data used in Appendix~\ref{app:eq}. Given a control polygon, we can determine whether a self-intersection of the B\'ezier curve occurs by determining whether $S$ (Equation~\ref{eq:S2}) has a root in $D$. We consider $S$ as a function $S(s,t,q)$ where $q=\{q_i\}_{i=0}^{n-1}$, so that finding a self-intersecting B\'ezier curve with an equilateral control polygon is equivalent to determining $s,t$ and $q$ such that the following are satisfied:

\begin{enumerate}
\item $S(s,t,q)=0$ where $s,t\in D$;
\item $||q_i||=r$ for each $i \in \{0,1,\cdots,n-1\}$, where $||\cdot||$ is the Euclidean norm;
\item $\sum_{i=0}^{n-1}{q_i}=\sum_{i=0}^{n-1}{P_{i+1}-P_i}=0$ since $P$ needs to be closed. 
\end{enumerate}

We assume $r=1$ without loss of generality since the value of $r$ can be adjusted by scaling. Throughout the provided codes, $n$ is always the degree of the B\'ezier curve. We give the code for function $S$ in Appendix~\ref{app:S}, where $[s,t]$ is labeled as $u$ and $q$ as $[a,b,c]$. 

The function $SF$ (Appendix~\ref{app:SF}) takes parameters for $s,t$ and $q$ and outputs a floating point value. It is designed to be zero if and only if the above three conditions are satisfied simultaneously. 

Precisely since $q$ should be taken from the unit sphere, $SF$ assigns $a,b,c$ values given by 
$$a=sin(\phi)cos(\theta); $$
$$b=sin(\phi)sin(\theta);$$
$$c=cos(\phi),$$
where $\phi$ and $\theta$ represent input parameters. In order to satisfy Condition $(3)$ above, $q_{n-1}$ is set equal to $-\sum_{i=0}^{n-2}{q_i}$. But in this way, $||q_{n-1}||$ may fail to be equal to $1$. So we include the function $F$ (Appendix~\ref{app:SF}) to determine whether $||q_{n-1}||=1$. The function is designed to be $F=||q_{n-1}||-1$ such that $||q_{n-1}||=1$ if and only if $F=0$. 

Symbolically, 
\begin{align}\label{eq:SF} SF=||S||+||F||,\end{align}
where $S$ is given by Equation~\ref{eq:S2} and $F=||q_{n-1}||-1$. Having the above three conditions satisfied simultaneously is equivalent to finding input values such that $SF=0$.

Since $SF \geq 0$, the minimum of $SF$ is $0$.  The function $SFminimizer$ (Appendix~\ref{app:SFminimizer}) uses $fminsearch$ (a numerical method integrated in MATLAB) to search for the minimum of $SF$, while returning this minimum and the corresponding values of $s,t$ and $q$. The user supplied initial guesses for $s,t$ and $q$ greatly influence the results, so we assign $SFminimizer$ randomized initial values $M$ times so that we get $M$ different minimums for different initial values. But no matter which initial values we use, as long as we can get ``a" minimum of $0$, then we get the equilateral control polygon which determines a self-intersecting B\'ezier curve. 

The data of finding the counterexample of Section~\ref{sec:ex} is included in Appendix~\ref{app:data}.

\section{Visualizing B\'ezier curves \& their control polygons}
\label{sec:vizb}

To study the knot types of a control polygon and its B\'ezier curve, a knot visualization tool was developed, called {\em Knot\_Spline\_Vis}.  The tool {\em Knot\_Spline\_Vis} takes a PL control polygon as input and displays both the PL curve and the associated B\'ezier curve.  For these studies, the input was always restricted to be a PL curve of known knot type.  The functions of {\em Knot\_Spline\_Vis} were designed to permit interactive studies of the topological relationships between a B\'ezier curve and its control polygon.  The intent is to use these examples to stimulate mathematical conjectures as a prelude to formal proofs.  

\begin{figure}[h!]
\centering
    \subfigure[Initial display]
    {
   \includegraphics[height=4cm]{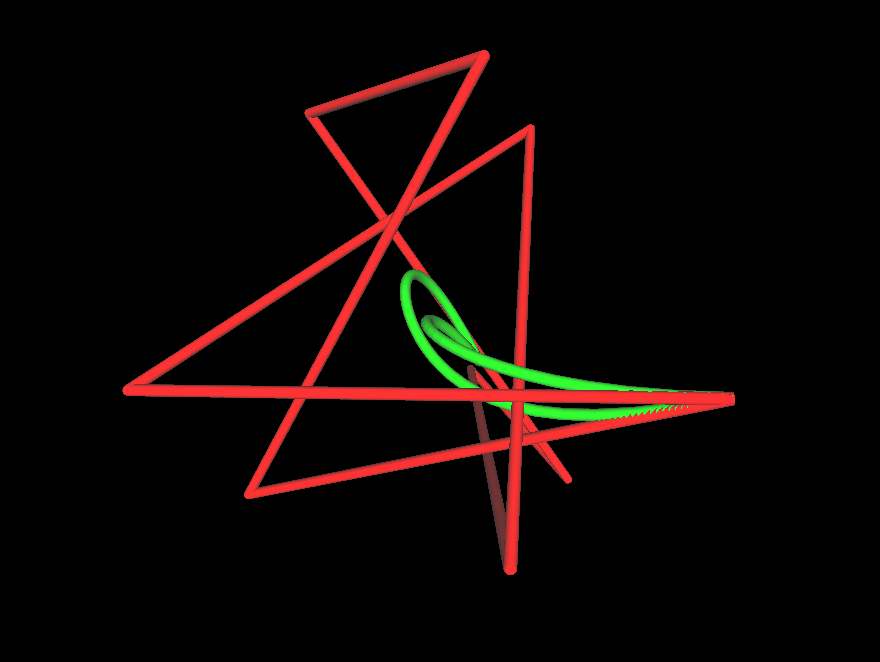}
    }
    \subfigure[Some vertices moved]
    {
   \includegraphics[height=4cm]{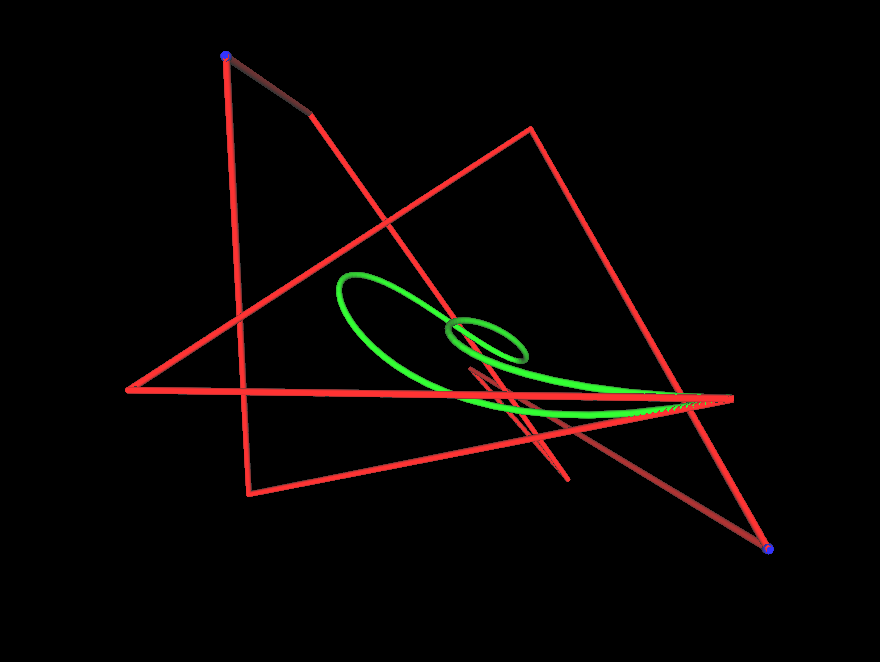}
    }

    \caption{Moving control points}
    \label{fig:move}
\end{figure}

Some of the standard graphical manipulation capabilities provided are illustrated in the following Figure~\ref{fig:rs} and~\ref{fig:zz}, where the rotation capabilities have been particularly useful to develop visual insights into the occurrence of self-intersections and crossings as fundamental for the study of knots. An editing window allows the user to change the coordinates of the control polygons, as shown in Figure~\ref{fig:move}. The software is freely available for download at the site \mbox{www.cse.uconn.edu/$\sim$tpeters}.  

\begin{figure}[h!]
\centering
   \includegraphics[height=5cm]{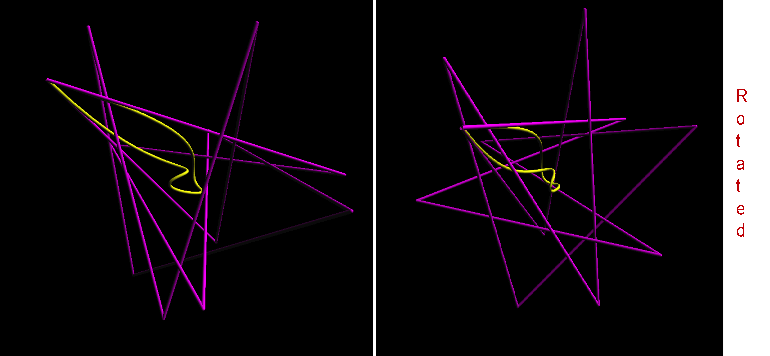}
    \caption{Rotating}
    \label{fig:rs}
\end{figure}

\begin{figure}[h!]
\centering
   \includegraphics[height=5cm]{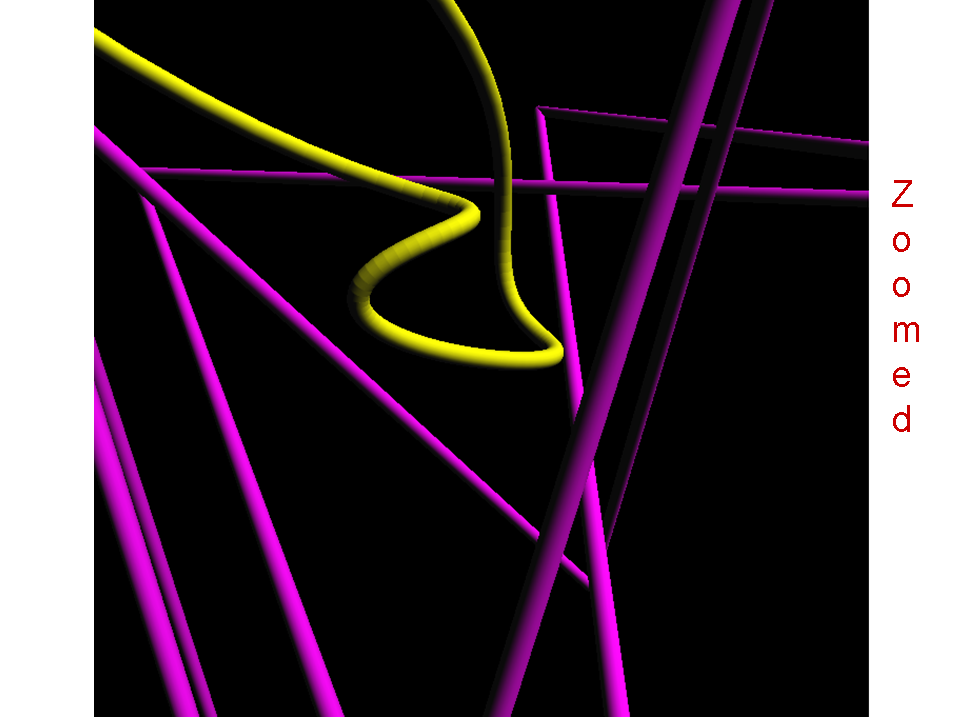}
    \caption{Scaling}
    \label{fig:zz}
\end{figure}

\begin{figure}[h!]
\centering
   \includegraphics[height=5cm]{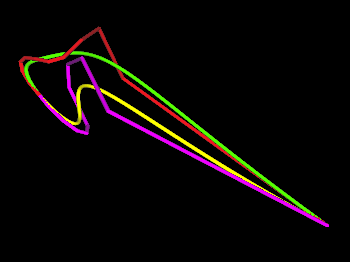}

    \caption{Subdivision}
    \label{fig:k10sub}
\end{figure}

The control polygon is an initial PL approximation to its associated B\'ezier curve.  Subdivision algorithms \cite{Piegl} produce additional control points to further refine the original control polygon, converging in distance to the B\'ezier curve.
Figure~\ref{fig:k10sub} illustrates performance of a standard subdivision technique, the de Casteljau algorithm \cite{G.Farin1990}. Users can specify a subdivision parameter. Figure~\ref{fig:k10sub} shows the case with parameter $\frac{1}{2}$.

\section{Conclusions and Future Work}
\label{sec:conc}

We use numerical algorithms and graphical software to generate counterexamples regarding the embeddings in $\mathbb{R}^3$ for a B\'ezier curve and its control polygon.  

First, we present a class of counterexamples showing an unknotted control polygon for a knotted B\'ezier curve and provide complete formal proofs of this condition.  We used a spline visualization tool, {\em Knot\_Spline\_Vis} for easy generation of many examples to gain intuitive understanding to formalize these proofs.  Secondly,  we used {\em Knot\_Spline\_Vis} in formulating the conjecture that any simple, closed equilateral control had a simple B\'ezier curve.  We provide contrary numerical evidence regarding this conjecture, as a useful guide for many applications in computer graphics and scientific visualization.

This work also shows the importance of continuing investigation 
\begin{itemize}
\item theoretically, into an infinite precision proof for a counterexample of the simple, closed equilateral polygon conjecture, where these demands will likely require further, substantive mathematical innovation and
\item practically, into a user interface for {\em Knot\_Spline\_Vis} to permit interactive 3D editing, as the current text driven interface was cumbersome.
\end{itemize}


\section{Web Posting of Supplemental Materials}
\label{sec:web}
Appendices listed below are posted on this webpage: 
\begin{center} \url{http://www.math.uconn.edu/~jili/Supplemental-materials.pdf} \end{center}

\noindent \textbf{Appendix A:} Numerical codes for knottedness of $\mathcal{C}$ (Figure~\ref{fig:k}):
\begin{itemize} \itemsep -2pt
\item A.1 Codes for searching self-intersections in Figure~\ref{fig:k-2d};
\item A.2 Codes for determining under or over crossings;
\item A.3 Data for searching intersections in Figure~\ref{fig:k-2d}.
\end{itemize}

\noindent \textbf{Appendix B:} Numerical codes for searching the Example of Figure~\ref{fig:eq}:
\begin{itemize} \itemsep -2pt
\item B.1 Codes for Equation~\ref{eq:S2};
\item B.2 Codes for Equation~\ref{eq:SF};
\item B.3 Codes for searching the roots of the system of Equations~\ref{eq:S2} and~\ref{eq:SF};
\item B.4 Data for finding the example shown in Figure~\ref{fig:eq}.
\end{itemize}

\pagebreak

\bibliographystyle{plain}
\bibliography{ji-tjp-biblio,tjp-add-biblio}

\begin{thebibliography}{10}

\bibitem{adams2004knot}
C.C. Adams.
\newblock {\em The Knot Book: An Elementary Introduction To The Mathematical
  Theory Of Knots}.
\newblock American Mathematical Society, 2004.

\bibitem{J.W.Alexander_G.B.Briggs1926}
J.~W. Alexander and G.~B. Briggs.
\newblock On types of knotted curves.
\newblock {\em Annals of Mathematics}, 28:562--586, 1926-1927.

\bibitem{Andersson1998}
L.~E. Andersson, T.~J. Peters, and N.~F. Stewart.
\newblock Selfintersection of composite curves and surfaces.
\newblock {\em CAGD}, 15:507--527, 1998.

\bibitem{Armstrong1983}
M.~A. Armstrong.
\newblock {\em Basic Topology}.
\newblock Springer, New York, 1983.

\bibitem{Bing1983}
R.~H. Bing.
\newblock {\em The Geometric Topology of 3-Manifolds}.
\newblock American Mathematical Society, Providence, RI, 1983.

\bibitem{Bisceglio}
J.~Bisceglio, T.~J. Peters, J.~A. Roulier, and C.~H. Sequin.
\newblock Unknots with highly knotted control polygons.
\newblock {\em CAGD}, 28(3):212--214, 2011.

\bibitem{G.Farin1990}
G.~Farin.
\newblock {\em Curves and Surfaces for Computer Aided Geometric Design}.
\newblock Academic Press, San Diego, CA, 1990.

\bibitem{Hass1999}
J.~Hass, J.~C. Lagarias, and N.~Pippenger.
\newblock The computational complexity of knot and link problems.
\newblock {\em Journal of the ACM}, 46(2):185--221, 1999.

\bibitem{JKSP10}
K.~E. Jordan, R.~M. Kirby, C.~Silva, and T.~J. Peters.
\newblock Through a new looking glass*: Mathematically precise visualization.
\newblock {\em SIAM News}, 43(5):1-- 3, 2010.

\bibitem{TCS08}
K.~E. Jordan, L.~E. Miller, E.~L.~F. Moore, T.~J. Peters, and A.~C. Russell.
\newblock Modeling time and topology for animation and visualization with
  examples on parametric geometry.
\newblock {\em Theoretical Computer Science}, 405:41--49, 2008.

\bibitem{JMPR11}
K.~E. Jordan, L.~E. Miller, T.~J. Peters, and A.~C. Russell.
\newblock Geometric topology and visualizing 1-manifolds.
\newblock In V.~Pascucci, X.~Tricoche, H.~Hagen, and J.~Tierny, editors, {\em
  Topology-based Methods in Visualization}, pages 1 -- 12, New York, 2011.
  Springer.

\bibitem{Livingston1993}
C.~Livingston.
\newblock {\em Knot Theory}, volume~24 of {\em Carus Mathematical Monographs}.
\newblock Mathematical Association of America, Washington, DC, 1993.

\bibitem{tjp-web}
T.~J. Peters and D.~Marsh.
\newblock Personal home page of {T. J. P}eters.
\newblock \url{http://www.cse.uconn.edu/}.

\bibitem{Piegl}
L.~Piegl and W.~Tiller.
\newblock {\em The NURBS Book}.
\newblock Springer, New York, 2nd edition, 1997.

\bibitem{KnotPlot}
R.~Scharein.
\newblock The knotplot site.
\newblock \url{http://www.knotplot.com/}.

\bibitem{Carlo}
C.~H. Sequin.
\newblock Spline knots and their control polygons with differing knottedness.
\newblock
  \url{http://www.eecs.berkeley.edu/Pubs/TechRpts/2009/EECS-2009-152.html}.

\end{thebibliography}

\pagebreak

\begin{appendices}

\section{Code: knottedness of  B\'ezier curve of Figure~\ref{fig:k}}

\subsection{Code for self-intersections of Figure~\ref{fig:k-2d}}\label{app:selfint}
\%The function $C2d(t)$ defines the projection of the B\'ezier curve. \\
function [value] = C2d(t)\\
n=10;\\
P=zeros(2,11);\\
P(1,:)= [-5.9, 10.3, -2.6, -10, 1.9, 11.2, -15.3, -11.7, 17.9, 2.9, -5.9];\\
P(2,:)= [4.7, -1.1, -12.4, 7, -12, 7.5, -1.7, 20, -1.1, -13.7, 4.7];\\
\\
sum=0;\\
for i=0:n\\
$sum=sum+nchoosek(n,i)* t^i * (1-t)^(n-i) * P(1,i+1)$;\\
end\\
$v_1=sum;$\\
\\
sum=0;\\
for i=0:n\\
$sum=sum+nchoosek(n,i)* t^i * (1-t)^(n-i) * P(2,i+1)$;\\
end\\
$v_2=sum$;\\
\\
$value = [v_1,v_2]$;\\
\\
\% Use the function `fminsearch' to find the minimums of $fnS(x)$, the zero minimums and parameters where the zeros occur are what we look for. \\
\\
function [value] = fnS(x)\\
$value = norm(C2d(x(1))-C2d(x(2)),2)$;\\
\\
function [value] = Smin(s0,t0)\\
Max=optimset('MaxFunEvals',1e+19);\\
comb=@(x)fnS(x);\\
u=[s0,t0];\\
$[xval,fval]=fminsearch(comb,u,Max)$

\subsection{Code for determining crossings}\label{app:cr}
\%Below is the function of the B\'ezier curve in 3D (Definition~\ref{def:c}). \\
function [value] = C(t)\\
n=10;\\
P=zeros(3,11);\\
P(1,:)= [-5.9, 10.3, -2.6, -10, 1.9, 11.2, -15.3, -11.7, 17.9, 2.9, -5.9];\\
P(2,:)= [4.7, -1.1, -12.4, 7, -12, 7.5, -1.7, 20, -1.1, -13.7, 4.7];\\
P(3,:)= [-6.2, 8.9, -6.3, -0.3, -0.6, -7.6, -4.1, 3.5, 2.9, 4.8, -6.2];\\
\\
sum=0;\\
for i=0:n\\
$sum=sum+nchoosek(n,i)* t^i * (1-t)^(n-i) * P(1,i+1);$\\
end\\
$v_1=sum$;\\
\\
sum=0;\\
for i=0:n\\
$sum=sum+nchoosek(n,i)* t^i * (1-t)^(n-i) * P(2,i+1);$\\
end\\
$v_2=sum$;\\
\\
sum=0;\\
for i=0:n\\
$sum=sum+nchoosek(n,i)* t^i * (1-t)^(n-i) * P(3,i+1);$\\
end\\
$v_3=sum$;\\
\\
$value = [v_1,v_2,v_3]$;

\subsection{Data for searching intersections in Figure~\ref{fig:k-2d}}\label{data:selfint}
\% The Matlab commands and corresponding results:\\
\% The initial input $(0.03,0.55)$ is figured out by the observation and calculations on self-intersections in Figure~\ref{fig:k-2d}. Similarly for the others below.\\
Smin(0.03,0.55)\\
xval = 0.0306    0.5573\\
fval =  2.9567e-04\\
\\
Smin(0.15,0.92)\\
xval = 0.1573    0.9244\\
fval = 1.5848e-04\\
\\
Smin(0.37,0.95)\\
xval = 0.3731    0.9493\\
fval = 1.4637e-04

\section{Code: generating counterexamples like Figure~\ref{fig:eq}}
\label{app:eq}

\subsection{Code for Equation~\ref{eq:S2}}\label{app:S}
function [value] = S(u,a,b,c)\\
    sum=0; subsum1=0; subsum2=0;\\
    for i=0:n-1\\
        for j=0:n-1-i\\
\% The use of adding $1$ to the vectors $a,b$ and $c$ is required by MATLAB.\\
            for k=0:i $subsum2=subsum2+nchoosek(i,k)*(1-u(2))^{(i-k)}*u(2)^k*a(j+k+1);$\\ 
            end\\
	   $subsum1=subsum1+subsum2 * nchoosek(n-1-i,j)*u(1)^{(n-1-i-j)}*(1-u(1))^j;$\\
        end\\
    $sum=sum+subsum1;$\\
    end\\
    $v_1=(1/n)*sum;$\\
\\
    sum=0; subsum1=0; subsum2=0;\\
    for i=0:n-1\\
        for j=0:n-1-i\\
            for k=0: i\\ 
            $subsum2=subsum2+nchoosek(i,k)*(1-u(2))^{(i-k)}*u(2)^k*b(j+k+1);$\\
	   end\\
	   $subsum1=subsum1+subsum2 * nchoosek(n-1-i,j)*u(1)^{(n-1-i-j)}*(1-u(1))^j;$\\
        end\\
    $sum=sum+subsum1;$\\
    end\\
    $v_2=(1/n)*sum;$\\
\\
    sum=0; subsum1=0; subsum2=0;\\
    for i=0:n-1\\
        for j=0:n-1-i\\
            for k=0:i $subsum2=subsum2+nchoosek(i,k)*(1-u(2))^{(i-k)}*u(2)^k*c(j+k+1);$\\
	   end\\
	    $subsum1=subsum1+subsum2 * nchoosek(n-1-i,j)*u(1)^{(n-1-i-j)}*(1-u(1))^j;$\\
        end\\
    $sum=sum+subsum1;$\\
    end\\
    $v_3=(1/n)*sum;$\\
$value = [v_1,v_2,v_3];$

\subsection{Code for Equation~\ref{eq:SF}}\label{app:SF}
\textbf{function [value] = SF(x,n)}\\
q=zeros(3,n); p=zeros(3,n+1); \\
a=zeros(1,n);b=zeros(1,n);c=zeros(1,n); \\
$alpha=zeros(1,2*n-2);$\\
for $i=1:2*n-2$\\
    alpha(i)=x(i);\\
end\\
for i=1:n-1\\
$a(i)=sin(alpha(i))*cos(alpha(n-1+i)); \\
b(i)=sin(alpha(i))*sin(alpha(n-1+i)); \\
c(i)=cos(alpha(i));$\\
q(:,i)=[a(i),b(i),c(i)];\\
end\\
\\
a(n)=0; b(n)=0; c(n)=0;\\
for i=1:n-1\\
a(n)=a(n)-a(i); b(n)=b(n)-b(i); c(n)=c(n)-c(i);\\
q(:,n)=[a(n),b(n),c(n)];\\
end\\
\\
u=zeros(1,2); $u=[x(2*n-1),x(2*n)];$\\
value=abs(F(alpha,n))+norm(S(u,a,b,c),2);\\
\\
for i=1:n\\
  for j=1:i;\\
      p(:,i+1)=p(:,i+1)+q(:,j);\\
  end\\
end\\
p\\
\\
\\
\textbf{function [value]=F(alpha,n)}\\
for i=1:n-1\\
$a(i)=sin(alpha(i))*cos(alpha(n-1+i)); b(i)=sin(alpha(i))*sin(alpha(n-1+i)); c(i)=cos(alpha(i));$\\
end\\
\\
a(n)=0; b(n)=0; c(n)=0;\\
for i=1:n-1\\
a(n)=a(n)-a(i); b(n)=b(n)-b(i); c(n)=c(n)-c(i);\\
end\\
$value=abs(a(n)^2+b(n)^2+c(n)^2-1);$

\subsection{Codes for solving the system of Equations~\ref{eq:S2} and~\ref{eq:SF}}
\label{app:SFminimizer}
function [value] = SFminimizer(n,M)\\
Max=optimset(\'MaxFunEvals\',1e+9);\\
comb=@(x)SF(x,n);\\
xval=zeros(2*n,M); \\
fval=zeros(1,M);\\
\\
for i=1:M\\  
phi=unifrnd(0,pi,1,n-1); theta=unifrnd(0,2*pi,1,n-1); \\
s0=unifrnd(0,1); t0=unifrnd(0,1-s0);\\
x0=[phi,theta,s0,t0];\\
$[xval(:,i),fval(i)]=fminsearch(comb,x0,Max);$\\
\%make sure $s+t=xval(2*n-1,i)+xval(2*n,i) < 1$ so that $s,t \in D$. \\
$if\ xval(2*n-1,i)+xval(2*n,i)> 1\ or\ xval(2*n-1,i)+xval(2*n,i)=1$\\
fval(i)=1;\\
else\\
end\\
end\\
\\
\%Display the minimums and corresponding $s,t$ and $q$. \\
xval fval\\
\%Returns the optimal one. \\
$[C,I]=min(fval);$\\
value=[C,I];

\subsection{Data for finding the example shown in Figure~\ref{fig:eq}}\label{app:data}
The counterexample (Section~\ref{sec:ex}) uses the second column below for the input parameters. \\
SFminimizer(6,10)\\
xval =\\
  Columns 1 through 10\\
    2.1907    0.0867    1.0279    2.4529    0.8294    1.3029    1.0958       0.8502    1.2136    0.3308\\
    2.1572    0.4353    1.9303    0.3841   -0.2788    0.3006   -0.2769      2.8797    3.0782    3.1699\\
    0.0079    3.2225    3.2107    1.8739    3.1482    3.6489    1.3900       0.5959    0.1228    1.8107\\
    0.5496    2.6633    0.2000    1.1955    0.5784    2.0316    1.6273       2.0583    3.1341    1.3276\\
    1.8511    2.9336    1.0424    2.6615    2.3439    1.2207    2.4580       4.3850    0.1877   -0.1146\\
    5.2200    1.2107    7.0598    3.4826    0.3696    0.5633    5.1432       4.2546    3.5523    5.5103\\
    2.7090    4.1128    4.1941    2.3445    2.4291    2.0589    0.4295       2.9467    1.0576    5.9281\\
    0.1954    2.0119    2.2644    7.0575    2.1168    2.3027    1.9744      -0.2264    2.5206    2.4902\\
    0.5312    7.4240    2.2923    5.1205    1.2132    4.0211    3.4493       1.6729    1.7930    5.7145\\
    1.3429    0.8465    1.4447    4.4661    4.0120    1.8860    0.0947       3.9222    5.1409    1.2114\\
    0.5200    0.2969    0.7942    0.7334    0.2149    0.6042    1.3269       0.9360    0.1185    0.0320\\
    0.0054    0.0633    0.2813    0.2543    0.6229    0.4588    0.0165       0.1667    0.5815    0.5970\\
\\
fval =\\
  Columns 1 through 7\\
    0.0000    0.0000    1.0000    0.2180    0.0001    1.0000    1.0000\\
  Columns 8 through 10\\
    1.0000    0.0001    0.0000\\
\\
\% The corresponding function values of $S$ and $F$ are given:\\
Svalue = $1.0e-03 *$\\
   -0.3861   -0.0970    0.1462\\
\\
Fvalue =$2.2329e-05$\\
\end{appendices}

\end{document}